\documentclass[10pt]{article}

\usepackage[latin1]{inputenc}
\usepackage{amsmath,amssymb}
\usepackage{amsmath,accents}
\usepackage{alltt,graphics,graphpap,color}
\usepackage[dvips]{graphicx}
\usepackage{amsfonts}
\usepackage{amscd}
\usepackage{mathrsfs}
\usepackage{hyperref}

\hypersetup{linktoc=page,colorlinks=true,linkcolor=red,  citecolor=blue, }

\newtheorem{Theorem}{Theorem}[section]

\newtheorem{Proposition}[Theorem]{Proposition}
\newtheorem{Lemma}[Theorem]{Lemma}
\newtheorem{Corollary}[Theorem]{Corollary}

\newtheorem{Notation}[Theorem]{Notation}

\newcommand{\ch}{\mathbb {CH} }
\newcommand{\hh}{\mathbb {HH} }
\newcommand{\sss}{\mathbb S}
\newcommand{\coq}{\cosh^2(\rho)}

\newcommand{\siq}{\sinh^2(\rho)}
\newcommand{\sih}{\sinh(\rho)}
\newcommand{\coh}{\cosh(\rho)}
\newcommand{\coqt}{\cosh^2(\rho(t))}
\newcommand{\siqt}{\sinh^2(\rho(t))}

\newcommand{\mm}{\mathcal {M}}

\newcommand{\gradphi}{\left|\nabla_{\sigma} \varphi\right|^2_{\sigma}}
\newcommand{\parho}{\frac{\partial}{\partial\rho}}
\newcommand{\pat}{\frac{\partial}{\partial t}}

\newcommand{\proof}{\noindent\emph{Proof. }}
\newcommand{\cvd}{\hfill$\square$ \bigskip}

\newcommand{\ep}{\varepsilon}

\newcommand{\tre}{\begin{array}{c|}3\\\hline\end{array}}
\begin{document}

%\def\qed{\hbox{\hskip 6pt\vrule width6pt height7pt
%depth1pt  \hskip1pt}\bigskip}

%\setlength\parindent{0pt}
%\vspace*{1 cm}

\title{Inverse mean curvature flow in quaternionic hyperbolic space }

\author{\sc Giuseppe Pipoli\footnote{Giuseppe Pipoli, Institut Fourier, Universit\'e Grenoble Alpes, 100 rue des Maths, 38610, Gi\` eres, France. E-mail: giuseppe.pipoli@univ-grenoble-alpes.fr}}
\date{}

\maketitle

%\vspace{3cm}
%\textbf{\Large 2st version}
%\vspace{1,5cm}

{\small \noindent {\bf Abstract:} In this paper we complete the study started in \cite{Pi2} of evolution by inverse mean curvature flow of star-shaped hypersurface in non-compact rank one symmetric spaces. We consider the evolution by inverse mean curvature flow of a closed, mean convex and star-shaped hypersurface in the quaternionic hyperbolic space. We prove that the flow is defined for any positive time, the evolving hypersurface stays star-shaped and mean convex. Moreover the induced metric converges, after rescaling, to a conformal multiple of the standard sub-Riemannian metric on the sphere defined on a codimension $3$ distribution. Finally we show that there exists a family of examples such that the qc-scalar curvature of this sub-Riemannian limit is not constant.} \medskip

\medskip

 \bigskip

\section{Introduction} \setcounter{equation}{0}

In this paper we complete the study started in \cite{Pi2} of evolution by inverse mean curvature flow of star-shaped hypersurface in non-compact rank one symmetric spaces, studying the evolution in the quaternionic hyperbolic space $\hh^n$. The inverse mean curvature flow is the most studied example in the class of expanding flows. For any given smooth hypersurface $F_0:\mm\rightarrow\hh^n$ , the solution of the inverse mean curvature flow with initial datum $F_0$ is the smooth one-parameter family of smooth immersions $F:\mm\times\left[0,T\right)\rightarrow\hh^n$ such that
\begin{equation}\label{imcf_def}
\left\{\begin{array}{rcl}
\displaystyle{\frac{\partial F}{\partial t}} & = & \displaystyle{\frac 1H \nu,}\\
F(\cdot,0) & = & F_0,
\end{array}\right.
\end{equation}
where $H$ is the mean curvature of $F_t=F(\cdot,t)$ and $\nu$ is the outward unit normal vector of $\mm_t=F_t(\mm)$. It is well known that this flow has a unique solution, defined at least for small times, if $\mm_0=F(\mm,0)$ is closed and mean convex. The behaviour of the evolution depends, of course, on the initial datum, but the geometry of the ambient manifold has a crucial role too. In fact Gerhardt \cite{Ge1} and and Urbas \cite{Ur} proved indipendently that, for any star-shaped hypersurface of the Euclidean space, the limit metric is, up to rescaling, always the standard round metric on the sphere. In \cite{HW} P.K. Hung and M.T. Wang showed that, when the ambient manifold is the hyperbolic space, the limit metric is not always round: it is a conformal multiple of the standard round metric on the sphere and so it is round only in special cases. More recently the author of the present paper studied the analogous problem in the complex hyperbolic space. In \cite{Pi2} it has been proved that in this case a new phenomenon appears: even after rescaling, the evolving metric blows up along a direction. Hence the limit metric is no more Riemannian, but only sub-Riemannian defined only on a codimension-$1$ distribution. Moreover there are infinite examples of initial data such that the limit sub-Riemannian metric does not have constant scalar Webster curvature. Recently similar problems have been studied in ambient manifolds which are warped products \cite{Sc,Zh}: $\mathbb{HH}^n$ is not in this class, in fact its metric can be written as \eqref{metrica}, but, for example, the Euclidean and the hyperbolic space are. Hence a different analysis is required. Other examples of evolution by inverse mean curvature flow of star-shaped hypersurfaces can be found in \cite{Di,Ne,KS}. The main result proved in this paper is the following.

\begin{Theorem}\label{main}
Let $n\geq 2$. For any $\mm_0$ closed, mean convex and star-shaped $\sss^3$-invariant hypersurface in $\hh^n$, let $\mm_t$ be its evolution by inverse mean curvature flow and let $g_t$ be the induced metric on $\mm_t$. Then:
\begin{itemize}
\item[(1)] $\mm_t$ is star-shaped for any time $t$;
\item[(2)] the flow is defined for any positive time;
\item[(3)] the rescaled induced metric $\tilde g_t=\left|\mm_t\right|^{-\frac {1}{2n+1}}g_t$ converges as $t\rightarrow\infty$ to $\tilde g_{\infty}=e^{2f}\sigma_{sR}$ for some smooth function $f$, where $\sigma_{sR}$ is the standard {sub-Riemannian} metric of the sphere $\sss^{4n-1}$ defined on a distribution of codimension $3$;
\item[(4)] finally there are examples of $\mm_0$ such that $\tilde g_{\infty}$ does not have constant qc-scalar curvature.
\end{itemize}
\end{Theorem}

The fact that the maximal time is infinite implies that the evolving hypersurface becomes arbitrary large and ``explores'' the structure at infinity of the ambient manifold as $t$ tends to infinity. Parts $(3)$ and $(4)$ of Theorem \ref{main} say that different initial data explore the structure at infinity of $\hh^n$ in different ways, but the conformal class of $\sigma_{sR}$ is preserved. 

Since we are considering only closed star-shaped hypersurfaces, we know that any hypersurface is an embedding of the sphere $\sss^{4n-1}$ into $\hh^n$. The rescaled induced metrics define a family of Riemannian metrics on $\sss^{4n-1}$. Analogously to what seen in \cite{Pi2}, this family of metrics diverges in some directions as the time goes to infinity. In our case we have $3$ independent special directions: $J_1\nu$, $J_2\nu$ and $J_3\nu$, where $J_1$, $J_2$ and $J_3$ are the complex structures of $\hh^n$ induced by the multiplication by the three quaternionic imaginary unities. Hence the limit metric is only sub-Riemannian defined on a distribution of codimension $3$. In this contest the best notion of curvature is the qc-curvature introduced by Biquard in \cite{Bi}. It is easy to prove that the evolution of a geodesic sphere stays a geodesic sphere at any time: it is a peculiar behaviour of the Euclidean space and rank one symmetric spaces. Then, in this special case, the qc-curvature of the limit rescaled metric is constant. On the other hand the research of a hypersurface with non trivial limit is more difficult. The qc-Yamabe problem was solved by Ivanov, Minchev and Vassilev in \cite{IMV}: they gave an explicit description of all the functions $f$ such that $e^{2f}\sigma_{sR}$ has constant qc-scalar curvature. With the help of this result we show that there are infinite examples of initial data such that $\tilde{g}_{\infty}$ has non constant qc-scalar curvature. Following the strategy introduced in \cite{Pi2}, the main tool is the following: for any star-shaped hypersurface $\mm$ we define the following Brown-York quantity
$$
Q(\mm)=|\mm|^{-1+\frac {1}{2n+1}}\int_{\mm}\left(H-\hat{H}\right)d\mu,
$$
where $|\mm|$ is the volume of $\mm$ and, if $\rho$ is the radial function defining $\mm$, $\hat{H}$ is the value of the mean curvature of a geodesic sphere of radius $\rho$ (see \eqref{hat_H} for the explicit definition). $Q$ gives a measure of how $\mm$ is far to being a geodesic sphere, however it isn't a true measure because, for example, $Q$ does not have a sign. In the final section of this paper we found the desired non-trivial example in the class of $\sss^3$-invariant submanifolds estimating the behaviour of $Q$ along the inverse mean curvature flow. 

This paper is organized as follows. In section 2 we collect some preliminaries and we fix some notation. In section 3 we compute the main geometric quantities for a star-shaped hypersurface in $\hh^n$, as the induced metric, the second fundamental form and the mean curvature. In section 4 we have ax explicit example, i.e. the evolution of the geodesic spheres. In section 5 we estimate the norm of the gradient of the radial function. As a consequence we have that the property of being star-shaped  and the mean convexity are preserved by the flow. In section 6 we collect without proof some results about the higher order derivatives of $\rho$  and we discuss their consequences. The proofs are not given because they are very similar of those of the analogous results in \cite{Pi2}. In particular we show that the solution of the flow is defined for any positive time and that the rescaled induced metric converges to a sub-Riemannian limit. Finally in section 7 we conclude the proof of Theorem \ref{main} studying the Webster curvature of the limit metric and giving a family of non-trivial examples.

%\begin{Remark}
%With the same techniques, the analogous of Theorem \ref{main} can be proved for closed, mean convex, star-shaped hypersurfaces in the octonion hyperbolic plane $\mathbb{OH}^2$: in this case we have a convergence to a sub-Riemannian metric defined on a distribution of codimension $7$. 
%\end{Remark}

\section{Preliminaries}  \setcounter{equation}{0}

\subsection{Riemannian and sub-Riemannian metrics on the sphere}
Every hypersurface considered in this paper is closed and star-shaped and so it is an embedding of $\sss^{4n-1}$, the sphere of dimension $4n-1$, into $\hh^n$. On that sphere we have different ``standard'' metrics. Let $\sigma$ be the usual round Riemannian metric on that sphere. We can distinguish three important vector fields: we can think the sphere embedded in $\mathbb H^n$ equipped with $J_1$, $J_2$ and $J_3$ the complex structures induced by the multiplication by the three quaternionic imaginary unities then:
$$
J_1^2=J_2^2=J_3^2=J_1J_2J_3=-id_{4n},
$$
where $id_{k}$ is the identity matrix of order $k$. If $\nu$ is the unit outward normal to $\sss^{4n-1}$, the $\xi_i=J_i\nu$ are three independent unit tangent vector field on the sphere. They are often called the \emph{Hopf vector fields}, because they span the tangent space of the fibers of the Hopf fibration : $\pi:\sss^{4n-1}\longrightarrow\mathbb{HP}^{n-1}$. They allow us to define the \emph{vertical distribution} $\mathcal V=\left\langle\xi_1,\xi_2,\xi_3\right\rangle$ and the \emph{horizontal distribution} $\mathcal H=\mathcal V^{\perp}$, the orthogonal complement of $\mathcal V$ with respect to $\sigma$. The \emph{Berger metrics} is a family of deformations of $\sigma$ obtained with the canonical deformation of $\pi$: for any $\lambda>0$ let $e_{\lambda}$ be the Riemannian metric defined by
$$
\left\{\begin{array}{rccl}
e_{\lambda}(X,Y) & = & \sigma(X,Y) &\quad\text{for any } X,\ Y\in \mathcal H;\\
e_{\lambda}(X,V) & = & 0 &\quad\text{for any } X\in \mathcal H,\ V\in\mathcal V;\\
e_{\lambda}(V,W) & = & \lambda\sigma(V,W)&\quad\text{for any } V,\ W\in \mathcal V.
\end{array}\right.
$$
When $\lambda$ converges to infinity, the metric $e_{\lambda}$ degenerates on the vertical directions, hence the limit is defined only on $\mathcal H$, but, since $\mathcal H$ is braket generating, $\sigma_{sR}=\lim_{\lambda\rightarrow\infty}e_{\lambda}$ is a sub-Riemannian metric on $\sss^{4n-1}$. We will call it the \emph{standard sub-Riemannian metric}. For brevity of notation, let us first define the following set: $\tre :=\left\{1,2,3\right\}$. The Levi-Civita connection of $e_{\lambda}$ has the following behaviour.

\begin{Lemma}\label{diff_LCC}
Fix a $\sigma$-orthonormal basis $(Y_1,\cdots,Y_{4n-1})$ of $\sss^{4n-1}$ such that for every $i\in\tre$ and for every $r$ $Y_i=\xi_i$  and  $Y_{4r+i}=J_iY_{4r}$. Let us denote with ${\nabla}_e$ ($\nabla_{\sigma}$ respectively) the Levi-Civita connection associated to the metric $e_{\lambda}$ ($\sigma$ respectively). Then for every $1\leq i,\ j\leq 2n-1$ we have:
$$
\nabla_{e\ Y_i}Y_j-\nabla_{\sigma\ Y_i}Y_j = \left\{\begin{array}{lll}
(1-\lambda)J_iY_j & \text{if }&i\in\tre,\ j\notin\tre; \\
 (1-\lambda)J_jY_i & \text{if }&j\in\tre,\ i\notin\tre;\\
0 & \text{otherwise.}
\end{array}\right.
$$
\end{Lemma}
\proof Fix $i\in\tre$ and $X$ any vector field tangent to $\sss^{4n-1}$, by the Gauss equation of the canonical immersion $\iota:\sss^{4n-1}\rightarrow\mathbb{R}^{4n}$ we have:
$$
\nabla_{\sigma X}J_i\nu=J_i\nabla_{0X}\nu-A_0(X,J_i\nu)\nu=-J_iX,
$$
where $\nabla_0$ is the Levi-Civita connection of the Euclidean space and $A_0$ is the second fundamental form of $\iota$. The thesis follows applying Lemma 3 of \cite{O} and Lemma 9.69 in \cite{Be}.
\cvd

\begin{Notation}\label{notazioni hess}
We introduce the following notation in order to distinguish between derivatives of a function with respect to different metrics. For any given function $f:\sss^{2n-1}\rightarrow\mathbb{R}$, let $f_{ij}$ ($\hat{f}_{ij}$ respectively) be the components of the Hessian of $f$ with respect to $\sigma$ ($e_{\ep}$ respectively). The value of $\ep$ will be clear from the context. The indices go up and down with the associated metric: for istance $\hat{f}_i^k=\hat{f}_{ij}e_{\ep}^{jk}$, while $f_i^k=f_{ij}\sigma^{jk}$. Analogous notations will be used for higher order derivatives.
\end{Notation}

\begin{Lemma}\label{hessiani}
Let $\varphi:\sss^{2n-1}\rightarrow\mathbb{R}$ be an $\sss^1$-invariant smooth function. With respect to the basis introduced in the previous Lemma, the Hessian of $\varphi$ with respect to $e_{\ep}$ is:
$$
\hat{\varphi}_{ij}=\left(\begin{array}{cc}
0 & \ep J_iY_j(\varphi)\\
\ep J_jY_i(\varphi)& {\varphi}_{ij}
\end{array}\right),
$$
where we are using Notations \ref{notazioni hess}. Taking the trace and the norm of the Hessian, in particular we have:
\begin{eqnarray*}
\Delta_e\varphi&=&\Delta_{\sigma}\varphi;\\
|\nabla^2_{e}\varphi|^2_{e} &=&|\nabla^2_{\sigma}\varphi|^2_{\sigma}+6(\ep-1)|\nabla_{\sigma}\varphi|^2_{\sigma}.
\end{eqnarray*}
\end{Lemma}
\proof
We can apply the previous result and we detect three cases:
\begin{itemize}
\item[1)] if $i,j\in\tre$
$$
\hat{\varphi}_{ij}=Y_iY_j(\varphi)-\nabla_{e Y_i}Y_j(\varphi)=0
$$
because $\varphi$ is $\sss^1$-invariant we have that $\varphi_j=Y_j(\varphi)=0$ for every $j\in\tre$ and 
$$
\nabla_{e Y_i}Y_j=\nabla_{\sigma Y_i}Y_j=\left\{\begin{array}{lll}
0 & \text{if} & i=j;\\
-J_iJ_j\nu & \text{if}& i\neq j.
\end{array}\right. 
$$
Hence $\nabla_{e Y_i}Y_j(\varphi)=0$ too.
\item[2)] if $i\in\tre$ and $j\notin\tre$ we get
$$
\hat{\varphi}_{ij}=\varphi_{ji}=Y_jY_i(\varphi)-\nabla_{e Y_j}Y_i(\varphi)=\mu J_iY_j.
$$
\item[3)] if $i,j\notin\tre$ it easy to check that $\hat{\varphi}_{ij}=\varphi_{ij}$.
\end{itemize}
Finally as a consequence of the symmetries considered, 
$$
|\nabla_e\varphi|^2_e=|\nabla_{\sigma}\varphi|^2_{\sigma}
$$
Taking into account this remark, the formulas for the Laplacian and the norm of the Hessian follow after some trivial computations.
\cvd

\subsection{The Biquard connection and the qc-scalar curvature}
The notion of \emph{quaternionic contact structure} (qc-structure for short) has been introduced by Biquard in \cite{Bi}. We refer also to the book \cite{IV1} of Ivanov and Vassilev for further details. A qc-structure on a real $(4n-1)$-dimensional manifold $\mm$ is a codimension $3$ distribution $\mathcal{H}$ (called \emph{horiziontal distribution}) locally given as the kernel of a $1$-form $\eta=(\eta_1,\eta_2,\eta_3)$ with values in $\mathbb{R}^3$ such that the three $2$-forms $d\eta_{i\left|\mathcal{H}\right.}$ are the fundamental forms of a quaternionic Hermitian structure on $\mathcal{H}$. Such $\eta$ is determined up to the action of $SO(3)$ on $\mathbb{R}^3$ and a conformal factor. Hence $\mathcal{H}$ is equipped with a conformal class $[g]$ of quaternionic Hermitian metrics. To every metric in the conformal class of $g$ one can associate a linear connection with torsion preserving the qc-structure called the \emph{Biquard connection}. It has been defined by Biquard in \cite{Bi} if $n>2$ and by Duchemin in \cite{Du} if $n=2$. They proved that this connection is unique. Using the Biquard connection we can define the qc-Ricci tensor as in \cite{Bi}: it is a symmetric tensor and its trace is called the \emph{qc-scalar curvature}. Once we know this new notion of curvature, we can define the \emph{qc-Yamabe problem}: determinate if there exists a metric $e^{2f}g$ in the conformal class $[g]$ with constant qc-scalar curvature. This problem has been solved in great generality, see \cite{IV2} for a survey. In our case $\mm=\sss^{4n-1}$, $g=\sigma_{sR}$ and for every $i$ $\eta_i(\cdot)=\sigma(\xi_i,\cdot)$. The metric $\sigma_{sR}$ has constant qc-curvature, but it is not the only one in its conformal class. In fact, in \cite{IMV}, Ivanov, Michev and Vassilev fully characterized the solution of the qc-Yamabe problem in the special case of the quaternionic Heisenberg group. As they noticed, with the Cayley transform, we can find the corresponding solutions on $\left(\sss^{4n-1},[\sigma_{sR}]\right)$: 
$e^{2f}\sigma_{sR}$ has constant qc-scalar curvature if and only if there are $c$ and $u$ positive constants and $\zeta\in\sss^{4n-1}$ such that
\begin{equation}\label{JL}
e^{-2f}(z)=c\left|\cosh(u)+\sinh(u)z\cdot\bar{\zeta}\right|^2, \quad\forall z\in\sss^{4n-1},
\end{equation}
Here we are considering the sphere of real codimension one immersed in $\mathbb{R}^{4n}\equiv\mathbb H^n$ and the norm and the product are the usual ones in $\mathbb H^n$.

The formula can be simplified imposing the symmetries considered. With some trivial computations similar to those of Lemma 2.5 of \cite{Pi2}, we can show that the following are equivalent:

\begin{Lemma}\label{caratterizzazione_curv}

Let $f:\sss^{4n-1}\rightarrow\mathbb R$ be an $\sss^3$-invariant function. 
\begin{itemize}
\item[(a)] $f$ satisfies \eqref{JL},
%\item[(b)] $e^{-f}$ is a linear combination of constants and first eigenfunctions of $\sss^{2n-1}$,
\item[(b)] $f$ is constant.
\end{itemize}
\end{Lemma}

\subsection{Quaternionic hyperbolic space}
The quaternionic hyperbolic space, like its real or complex analogous, can be defined in many equivalent ways. Since we wish to study star-shaped hypersurfaces, it is convenient introduce in polar coordinates. Let $\hh^n$ be $\mathbb R^{4n}$ equipped with the following $\bar g$: 
\begin{equation}\label{metrica}
\bar g=d\rho^2+\sinh^2(\rho)e_{\cosh^2(\rho)},
\end{equation}

where $\rho$ represents the distance from the origin and $e_{\cosh^2(\rho)}$ is the Berger metric of parameter $\coq$ on $\sss^{4n-1}$.  Note that $\bar g$ is not a warped product, so this ambient space is not included in the case studied in \cite{Sc, Zh}. Since the inverse mean curvature flow is invariant under the action of an isometry of the ambient space, all the hypersurfaces considered in this paper can be thought as star-shaped with respect to the origin of the polar coordinates. The metric $\bar g$ has some nice properties: it is an example of rank one symmetric space, hence $\bar \nabla\bar R=0$. 

Here and in the following we are using the convention to put a bar over the symbol for geometric quantity of the fixed ambient manifold $\hh^n$, for example $\bar{\nabla}$ is the Levi-Civita connection of $\bar{g}$.

Moreover its Riemann curvature tensor has the following explicit expression
\begin{equation}\label{curv}
\begin{array}{rcl}
\bar R(X,Y,Z,W)&=&-\bar g(X,Z)\bar g(Y,W)+\bar g(X,W)\bar g(Y,Z)\\&&-\sum_{i=1}^3\bar g(X,J_iZ)\bar g(Y,J_iW)+\bar g(X,J_iW)\bar g(Y,J_iZ)\\
&&-2\sum_{i=1}^3\bar g(X,J_iY)\bar g(Z,J_iW),
\end{array}
\end{equation}
where $J_1$, $J_2$ and $J_3$ are the complex structure of $\hh^n$ induced by the quaternionic imaginary units. It follows that the sectional curvature of a plane spanned by two orthonormal vectors $X$ and $Y$ is given by
\begin{equation}\label{sez}
\bar K(X\wedge Y) = -1-3\sum_{i=1}^3\bar g(X,J_iY)^2=-1-3|pr_{Y\mathbb{H}}X|^2,
\end{equation}
where $Y\mathbb{H}$ is the space spanned by $J_iY$ and $pr_{Y\mathbb{H}}$ is the projection on that space.
Then $-4\leq\bar K\leq -1$ and it is equal to $-1$ (respectively to $-4$) if and only if $X$ is orthogonal (respectively parallel) to $Y\mathbb{H}$. Furthermore $\hh^n$ is Einstein with $\bar Ric=-4(n+2)\bar g$.

\subsection{Inverse mean curvature flow}
Since we are considering only closed and mean convex initial data, the standard theory ensures that the inverse mean curvature flow \eqref{imcf_def} has a unique smooth solution, at least for small times. Here we list the evolutions of the main geometric quantities. The proof of this Lemma is similar to the computation of the analogous equations for the mean curvature flow which can be found in \cite{Hu}. We use the following notations: let $g_{ij}$ be the induced metric, and $g^{ij}$ its inverse; the second fundamental form is denoted with $h_{ij}$, while the mean curvature is $H=h_{ij}g^{ji}$. Finally $\left|\mm_t\right|$ denotes the volume.

\begin{Lemma}\label{eq_evoluz} Since the ambient space is symmetric the following evolution equations hold:
\begin{itemize}
\item[(1)] $\displaystyle{\frac{\partial g_{ij}}{\partial t} = \frac 2H h_{ij}}$,
\item[(2)] $\displaystyle{\frac{\partial g^{ij}}{\partial t} = -\frac 2H h^{ij}}$,
\item[(3)] $\displaystyle{\frac{\partial H}{\partial t} = \frac{\Delta H}{H^2}-2\frac{\left|\nabla H\right|^2}{H^3}-\frac{\left|A\right|^2}{H}-\frac{\bar{R}ic(\nu,\nu)}{H}}$,
%\item[(4)] $\displaystyle{\frac{\partial h_{ij}}{\partial t} =\frac{\nabla_i\nabla_jH}{H^2}+\frac{h_i^lh_{lj}}{H}-\frac{2}{H^3}\nabla_iH\nabla_jH-\frac{1}{H}\bar R_{i0j0}}$,
\item[(4)] $\displaystyle{\frac{\partial h_{ij}}{\partial t} =\frac{\Delta h_{ij}}{H^2}-\frac{2}{H^3}\nabla_iH\nabla_jH+\frac{\left|A\right|^2}{H^2}h_{ij}-\frac{2}{H}\bar R_{i0j0}+\bar{R}ic(\nu,\nu)\frac{h_{ij}}{H^2}}$
\item[\phantom{(4')}]$\displaystyle{\phantom{\frac{\partial h_{ij}}{\partial t} =}+\frac{1}{H^2}g^{lr}g^{ms}\left(2\bar R_{risj}h_{lm}-\bar R_{rmis}h_{jl}-\bar R_{rmjs}h_{il}\right)}$,
\item[(5)] $\displaystyle{\frac{\partial h_i^j}{\partial t} =\frac{\Delta h_i^j}{H^2}-\frac{2}{H^3}\nabla_iH\nabla_kHg^{kj}+\frac{\left|A\right|^2}{H^2}h_i^j-\frac{2}{H}\bar R_{i0k0}g^{kj}-2\frac{h_i^kh_k^j}{H}}$
\item[\phantom{(5)}]$\displaystyle{\phantom{\frac{\partial h_i^j}{\partial t} =}+\bar{R}ic(\nu,\nu)\frac{h_i^j}{H^2}+\frac{1}{H^2}g^{lr}g^{ms}g^{kj}\left(2\bar R_{risk}h_{lm}-\bar R_{rmis}h_{kl}-\bar R_{rmks}h_{il}\right)}$,
\item[(6)] $\displaystyle{\frac{d\left|\mm_t\right|}{d t}= \left|\mm_t\right|}$,
\item[(7)] $\displaystyle{\frac{\partial\nu}{\partial t}= \frac{\nabla H}{H^2}}$.
\end{itemize}
\end{Lemma}

Here and in the following we are using Einstein convention on repeted indices. Moreover the operation of raising/lowering the indices is done with respect to the metric: for example $h_i^j=h_{ik}g^{kj}$.  Note that, integrating equation $(6)$, we have that the inverse mean curvature flow is en expanding flow, precisely $\left|\mm_t\right|=\left|\mm_0\right|e^t$. 

Moreover, with the same proof of Lemma 3.1 of \cite{Pi1} we can prove the following result. 

\begin{Lemma}
The evolution of an $\sss^3$-invariant hypersurface stays $\sss^3$-invariant during the whole duration of the flow.
\end{Lemma}
\section{Geometry of star-shaped hypersurfaces}\setcounter{equation}{0}

In this section we compute the main geometric quantities for a star-shaped hypersurface in $\hh^n$. 
Let $F:\sss^{4n-1}\to\hh^n$ be a smooth star-shaped immersion. Up to an isometry of the ambient space, we can consider that it is star-shaped with respect to the origin. Then $F$ is defined by its radial function: there exist a smooth function $\rho:\sss^{4n-1}\to\mathbb R^+$ such that in polar coordinates $\mm=F(\sss^{4n-1})=\left\{(z,\rho(z))\in\hh^n\left|z\in\sss^{4n-1}\right.\right\}$. Fix any $(Y_1,\dots,Y_{4n-1})$ tangent basis of the sphere $\sss^{4n-1}$, for every $i$ we define $\rho_i=Y_i(\rho)$ and $V_i=F_*Y_i\equiv Y_i+\rho_i\frac{\partial}{\partial \rho}$. Then $(V_1,\dots,V_{4n-1})$ is a tangent basis of $\mm$. The induced metric on $\mm$ is $g=F^*\bar g$, in local coordinates we have
\begin{equation}\label{metrica}
g_{ij}=\rho_i\rho_j+\siq e_{ij},
\end{equation}
where for short $e_{ij}=(e_{\coq})_{ij}$. The inverse of the metric therefore is
\begin{equation}\label{metrica_inv}
g^{ij}=\frac{1}{\siq}\left(e^{ij}-\frac{\rho^i\rho^j}{\siq+|\nabla_e\rho|^2_e}\right),
\end{equation}
where $e^{ij}$ is the inverse of $e_{ij}$, $\rho^i=\rho_ke^{ki}$ and the gradient and the norm of $\rho$ are defined with respect to the metric $e_{\coq}$. In order to simplify the expressions we can fix a function $\varphi=\varphi(\rho)$ such that $\frac{d\varphi}{d\rho}=\frac{1}{\sinh(\rho)}$ and introduce $v=\sqrt{1+\frac{|\nabla_e\rho|^2_e}{\siq}}$. Since $\varphi_i=Y_i(\varphi)=\frac{\rho_i}{\sinh(\rho)}$, we get  $g_{ij}=\siq(\varphi_i\varphi_j+e_{ij})$, $g^{ij}=\frac{1}{\siq}\left(e^{ij}-\frac{\varphi^i\varphi^j}{v^2}\right)$ and $v=\sqrt{1+|\nabla_e\varphi|^2_e}$. A unit normal vector is $\nu=\frac{1}{v}\left(\parho-\frac{\nabla_e\rho}{\siq}\right)$. In case of $\sss^3$-invariant hypersurfaces we have that $\nabla_e\rho=\nabla_{\sigma}\rho$ and $\nabla_e\varphi=\nabla_{\sigma}\varphi$.

Now we want to compute the second fundamental form of $\mm$. Since the metric of the ambient space is not the same in any direction, it is convenient to introduce a specific basis tangent to $\sss^{4n-1}$.  Then let us use the basis $(Y_1,\dots,Y_{4n-1})$ introduced in Lemma \ref{diff_LCC}. In this way we have

\begin{equation}\label{berger}
e_{\coq}=\left(\begin{array}{cc}
\coq id_3& 0 \\
0 & id_{4n-4}
\end{array}
\right),
\end{equation}
For each $i$ and $j$ let $h_{ij}=-\bar g\left(\bar\nabla_{V_i}V_j,\nu\right)$. Moreover we introduce the following notation: Latin indices $i,j,...$ range from $1$ to $4n-1$ and are related to components tangent to the sphere, the index $0$ represents the radial direction $\parho$ and Greek indices $\alpha, \beta,...$ range from $0$ to $4n-1$. An explicit computation, together to the fact that for every $i$ $\parho(\rho_i)=0$ and $\bar g(\bar\nabla_{Y_i}\parho,\parho)=0$ we get:

$$
h_{ij}=\frac 1v\left(\bar{\Gamma}_{ij}^k\rho_k+\rho_i\rho_k\bar{\Gamma}_{0j}^k+\rho_j\rho_k\bar{\Gamma}_{0i}^k-\bar{\Gamma}_{ij}^0-Y_i(\rho_j)\right).
$$
We have that $\bar{\Gamma}_{ij}^k=\hat{\Gamma}_{ij}^k$, the Christoffel symbols of the metric $e_{\coq}$, then $\bar{\Gamma}_{ij}^k\rho_k-Y_i(\rho_j)=-\hat{\rho}_{ij}$, where we used Notation \ref{notazioni hess}. For short, let $Y_0$ be $\parho$, then

\begin{eqnarray*}
\bar{\Gamma}_{ij}^0&=&\frac 12\left(Y_i(\bar g_{i\alpha})-Y_{\alpha}(\bar g_{ij})+Y_j(\bar g_{i\alpha})\right)\bar g^{\alpha 0}\\
&= & -\frac 12\parho\left(\bar g_{ij}\right)=-\gamma_i\delta_{ij}=\left\{\begin{array}{ll}
-\sih\coh(\siq+\coq)\delta_{ij} & \text{if } i\in\tre\ ,\\
-\sih\coh\delta_{ij} & \text{if } i\notin\tre\ .  
\end{array}\right.
\end{eqnarray*}

\noindent Finally
\begin{eqnarray*}
\bar{\Gamma}_{i0}^k & = & \frac 12 \left(Y_i\left(\bar g_{0\alpha}\right)-Y_{\alpha}\left(\bar g_{i0}\right)+\parho\left(\bar g_{i\alpha}\right)\right)\bar g^{\alpha k}\\
& =& \frac 12 \parho\left(\bar g_{i\alpha}\right)\bar g^{\alpha k}=\eta_i\delta_{ik}
=\left\{	\begin{array}{ll}
\displaystyle{\frac{\siq+\coq}{\sinh(\rho)\cosh(\rho)}\delta_{ik} }& \text{if } i\in\tre\ ,\\
\displaystyle{\frac{\cosh(\rho)}{\sinh(\rho)}\delta_{ik} }& \text{if } i\notin\tre\ .  
\end{array}\right.
\end{eqnarray*}

Note that $$\hat{\varphi}_{ij}=\frac{1}{\sinh(\rho)}\rho_{ij}-\frac{\cosh(\rho)}{\siq}\rho_i\rho_j\quad\Leftrightarrow\quad\hat{\rho}_{ij}=\sinh(\rho)\varphi_{ij}+\sinh(\rho)\cosh(\rho)\varphi_i\varphi_j.$$
Similar equations hold for $\rho_{ij}$ and $\varphi_{ij}$ too.

Summing up these quantities we get

%\begin{equation}\label{sff}
%h_{ij}  = -\frac{\sih}{v}\varphi_{ij} +\frac{\sih}{v}*
%\left\{\begin{array}{ll}
%\frac{2\siq+\coq}{\coh}\varphi_i\varphi_j+\coh(\siq+\coq)\delta_{ij} & \text{ if } i,j\in\tre\ ,\\
%\coh\varphi_i\varphi_j+\coh\delta_{ij} & \text{ if } i,j\notin\tre\ ,\\
%\frac{\siq+\coq}{\coh}\varphi_i\varphi_j+\coh(\siq+\coq)\delta_{ij} & \text{ otherwise.}
%\end{array}
%\right.
%\end{equation}
%\textbf{delta o e?}

\begin{equation}\label{sff}
h_{ij}  = -\frac{\sih}{v}\hat{\varphi}_{ij} +
\left\{\begin{array}{ll}
\displaystyle{\frac{\siq+\coq}{v\sih\coh}g_{ij}+\frac{\sinh^3(\rho)}{v\coh}\varphi_i\varphi_j }& \text{ if } i,j\in\tre\ ,\\
\displaystyle{\frac{\coh}{v\sih}g_{ij} }& \text{ if } i,j\notin\tre\ ,\\
\displaystyle{\frac{\siq+\coq}{v\sih\coh}g_{ij}}& \text{ otherwise.}
\end{array}
\right.
\end{equation}

Raising the second index we have
\begin{equation}\label{sff_up_1}
h_{i}^k  = -\frac{\hat{\varphi}_{ij}\tilde{e}^{jk}}{v\sih} +
\left\{\begin{array}{ll}
\displaystyle{\frac{\siq+\coq}{v\sih\coh}\delta_i^k+\frac{\sih}{v\coh}\sum_{j=1}^3\varphi_i\varphi_j\tilde{e}^{jk}}& \text{ if } i\in\tre\ ,\\
\displaystyle{\frac{\coh}{v\sih}\delta_{i}^k +\frac{\sih}{v\coh}}\sum_{j=1}^3\tilde{e}_{ij}\tilde{e}^{jk}& \text{ if } i\notin\tre\ ,
\end{array}
\right.
\end{equation}
where $\tilde e^{ij}=\siq g^{ij}=e^{ij}-\frac{\varphi^i\varphi^j}{v^2}$
Taking the trace of \eqref{sff_up} we obtain the mean curvature of $\mm$:

\begin{eqnarray}
\nonumber H & = & h_i^i\\
\nonumber &=&-\frac{\hat{\varphi}_{ij}\tilde e^{ji}}{v\sih}+3\frac{\siq+\coq}{v\sih\coh}+4(n-1)\frac{\coh}{v\sih}\\
\nonumber && +\frac{\sih}{v\coh}\sum_{k=1}^3\varphi_{i}\varphi_k\tilde{e}^{ik}+\frac{\sih}{v\coh}\sum_{k=1}^3\tilde{e}_{ik}\tilde{e}^{ik}\\
%\nonumber & =& -\frac{\varphi_{ij}\tilde e^{ji}}{v\sih}+3\frac{\sih}{v\coh}+(4n-1)\frac{\coh}{v\sih}\\
%\nonumber && +3\frac{\sih}{v\coh}\tilde{e}_{i1}\tilde{e}^{i1}-\frac{\sih\coh}{v}\sum_{k=1}^3\tilde{e}^{kk}\\
& = & -\frac{\hat{\varphi}_{ij}\tilde e^{ji}}{v\sih}+\frac{\hat{H}}{v}+\frac{\sih}{v^3\coh}\sum_{k=1}^3(\varphi_k)^2,\label{mc_1}
\end{eqnarray}
where \begin{equation}\label{hat_H}
\hat{H}=\hat{H}(\rho)=(4n-1)\frac{\coh}{\sih}+3\frac{\sih}{\coh}%=\frac{2n\coq-1}{\sih\coh}.
\end{equation}
For the sake of clarity, we point out that in the above formula, the Einstein convention is used for all repeated indices except that, for $k$ which belongs to $\tre$. 

In case of $\sss^3$-invariance the formula just found can be simplified:
\begin{equation}\label{sff_up}
h_{i}^k  = -\frac{\hat{\varphi}_{ij}\tilde{e}^{jk}}{v\sih} +
\left\{\begin{array}{ll}
\displaystyle{\frac{\siq+\coq}{v\sih\coh}\delta_i^k}& \text{ if } i\in\tre\ ,\\
\displaystyle{\frac{\coh}{v\sih}\delta_{i}^k}& \text{ if } i\notin\tre.
\end{array}
\right.
\end{equation}
Moreover let $\tilde{\sigma}^{ji}=\sigma^{ji}-\frac{\varphi^i\varphi^j}{v^2}$, then by Lemma \ref{hessiani} we have that $\hat{\varphi}_{ij}\tilde{e}^{ji}=\varphi_{ij}\tilde{\sigma}^{ji}$. It follows that:
\begin{equation}
H= -\frac{\varphi_{ij}\tilde {\sigma}^{ji}}{v\sih}+\frac{\hat{H}}{v}.\label{mc}
\end{equation}

\section{The case of geodesic spheres}

If the radial function $\rho$ is constant, $\mm_0$ is, of course, a geodesic sphere. From \eqref{sff_up} we can see that $\mm_0$ has two distinct principal curvatures: $\lambda(\rho)=\frac{\coh}{\sih}$ with multiplicity $4n-4$ and $\mu(\rho)=2\coth(2\rho)=\frac{\sih}{\coh}$ with multiplicity $3$ and eigenvectors $\xi_i$, with $i\in\tre$. It follows that the mean curvature is costant and depends only on the radius. Then the evolution of $\mm_0$ by inverse mean curvature flow reduces to an ODE: the evolution on $\mm_0$ is a family of geodesic spheres $\mm_t$ of radius $\rho(t)$ satisfying
$$\left\{\begin{array}{l}
\dot{\rho}=\frac 1H=\frac{\sih\coh}{(4n-1)\coq+3\siq},\\
\rho(0) = \rho_0.
\end{array}\right.
$$
Trying to solve this ODE we can see that the solution is defined for any positive time and $\rho(t)=\frac{t}{2(2n+1)}+o(1)$ as $t\rightarrow\infty$. The rescaled induced metric is
$$
\tilde g_t=\left|\mm_t\right|^{-\frac{1}{2n+1}}g=\frac{\siqt}{\left|\mm_0\right|^{\frac {1}{2n+1}} e^{\frac {t}{2n+1}}}e_{\coqt}.
$$
Since in this special case $\rho$ depends only on time, we have that $\tilde{g}_t$ converges to $\tilde{g}_{\infty}=c^2\sigma_{sR}$, where $c$ in a constant different from zero. We can see for the first time the phenomenon announced in the statement of Theorem \ref{main}: the rescaled metric converge to a sub-Riemannian metric defined only on a distribution of codimension $3$. Moreover the qc-curvature of this metric is constant.

The following result is useful to bound the evolution of the radial function in the general case.

\begin{Lemma}\label{confronto}
Consider two concentric geodesic spheres in $\ch^n$ of radius $\rho_1(0)$ and $\rho_2(0)$ respectively. For every $i=1,2$, let $\rho_i(t)$ be the evolution by inverse mean curvature flow of initial datum $\rho_i(0)$, then there exist a positive constant $c$ depending only on $n$, $\rho_1(0)$ and $\rho_2(0)$ such that for every time we have
$$
\left|\rho_2(t)-\rho_1(t)\right|\leq c\left|\rho_2(0)-\rho_1(0)\right|.
$$
\end{Lemma}

\proof Suppose that $\rho_2(0)>\rho_1(0)$ and let us define $\delta=\delta(t)=\rho_2(t)-\rho_1(t)$. The proof is similar to that of Lemma 4.1 of \cite{Pi2}, taking into account that this time the function $\delta$ satisfies
\begin{eqnarray*}
\frac{d\delta}{dt} & = & \frac{1}{(4n-1)\coth(\rho_2)+3\tanh(\rho_2)}-\frac{1}{(4n-1)\coth(\rho_1)+3\tanh(\rho_1)}.
\end{eqnarray*}
\cvd

From these properties of the geodesic spheres we can deduce some estimates on the evolution of the general case: we can bound the inverse mean curvature flow of a general star-shaped hypersurface $\mm_0$ with the evolution of geodesic sphere inside and a geodesic sphere outside. Applying Lemma \ref{confronto} we have that the oscillation of the radial function of $\mm_0$ is bounded by a constant which depends only on the initial datum. Below we will show that the flow is defined for any positive time also for any star-shaped initial datum. It follows that in any case we have $\rho(t)=\frac{t}{2(2n+1)}+o(1)$ as $t\rightarrow\infty$.

\section{First order estimates}\setcounter{equation}{0}
The main result of this section is the proof of part $(1)$ of Theorem \ref{main}. Moreover we will prove also that the mean curvature converges exponentially to that one of an horosphere. The main technical result is the following:

\begin{Proposition}\label{gradphi}
There exist a positive constants $c$ such that
$$
\gradphi\leq ce^{- \frac{t}{2n+1}}.
$$

\end{Proposition}

As an immediate geometric consequence we have:

\begin{Corollary}\label{stellato}
The evolution of any star-shaped hypersurface stays star-shaped for any time the flow is defined.
\end{Corollary}
\proof An hypersurface is star-shaped if and only if $\parho$ and $\nu$ are never orthogonal in $\mathbb{HH}^n$. This means that there exists a positive constant $c$ such that
$$
\bar g\left(\parho,\nu\right)=\frac 1v\geq \frac 1c \quad\Leftrightarrow\quad v^2\leq c^2
$$
Recalling that $v^2=1+\gradphi$, the thesis follows from Proposition \ref{gradphi}.
\cvd

The proof of Theorem \ref{gradphi} is divided into two steps: first we can prove that $\gradphi$ is bounded, then, applying the strategy of Section 5 of \cite{Pi2}, we can prove the exponential decay.

\begin{Lemma}\label{gradphi_w}
$$\gradphi\leq\sup_{z\in\sss^{4n-1}}\left|\nabla_{\sigma}\varphi(z,0)\right|^2_{\sigma}.$$
\end{Lemma}
\proof Let us define $\omega=\frac 12 \gradphi$. We want to compute the evolution equation for $\omega$.
We have
%$$
%\frac {1}{Hv}=\frac{d\rho}{dt}=\frac{\partial\rho}{\partial t}+\frac{\partial\rho}{\partial x_i}\frac{\partial x_i}{\partial t}=\frac{\partial\rho}{\partial t}-\frac{\rho^i\rho_i}{Hv\siq}.
%$$
%Then
\begin{equation}\label{ev_rho}
\frac{\partial\rho}{\partial t}=\frac vH,
\end{equation}
and so
\begin{equation}\label{ev_phi}
\frac{\partial\varphi}{\partial t}=\frac {1}{\sih}\frac{\partial\rho}{\partial t}=\frac{v}{H\sih}=:\frac 1 F
\end{equation}
holds. From the explicit computation of the mean curvature \eqref{mc} we have
$$
F=F(\varphi,\varphi_i,\varphi_{ij})=\frac{-\varphi_{ij}\tilde {\sigma}^{ij}}{v^2}+\frac{\sih\hat{H}}{v^2}.
$$
Now we can compute the evolution equation of $\omega$: let $a^{ij}=-\frac{\partial F}{\partial\varphi_{ij}}=\frac{\tilde {\sigma}^{ij}}{v^2}$ and $a^i=\frac{\partial F}{\partial\varphi_{i}}$, then
\begin{eqnarray*}
\frac{\partial\omega}{\partial t}&=&\varphi^k \nabla_k\frac{\partial\varphi}{\partial t}\\
& =& -\frac{1}{F^2}\left( -a^{ij}\varphi_{ijk}\varphi^k+a^i\varphi_{ik}\varphi^k+\frac{\partial F}{\partial\varphi}\varphi_k\varphi^k\right)\\
& = &-\frac{1}{F^2}\left( -a^{ij}\varphi_{ijk}\varphi^k +a^i\omega_i+2\frac{\partial F}{\partial\varphi}\omega\right)
\end{eqnarray*}
We can apply the rule interchanging for derivatives:
$$
\varphi_{ijk}=\varphi_{kji}+R^m_{\phantom{m}ijk}\varphi_m,
$$
where this time $R$ is the Riemannian curvature tensor of $\sigma$, i.e. $R_{sijk}=\sigma_{sj}\sigma_{ik}-\sigma_{sk}\sigma_{ij}$. Since $a^{ij}$ is symmetric we get:
\begin{eqnarray*}
 -a^{ij}\varphi_{ijk}\varphi^k &= & -a^{ij}\varphi_{kji} \varphi^k-a^{ij}\left(\delta_j^m\sigma_{ik}-\delta_k^m\sigma_{ij}\right)\varphi_m \varphi^k\\
& = & -a^{ij}\omega_{ij}+a^{ij}\varphi_{ki} \varphi_j^k-a^{ij}\varphi_i \varphi_j+2a^i_i\omega\\
&\geq &-a^{ij}\omega_{ij}.
\end{eqnarray*}
Moreover we have
\begin{eqnarray}
\nonumber\frac{\partial F}{\partial\varphi} & = & \frac{\partial F}{\partial\rho}\frac{\partial \rho}{\partial\varphi}=\frac{\sih}{v^2}\frac{\partial}{\partial\rho}\left(\sih\hat{H}\right)\\
  & = & \frac{\siq}{v^2\coq}\left(2(2n+1)\coq+3\right).\label{stima01}
\end{eqnarray}
In particular $\frac{\partial F}{\partial\varphi}>0$. We have that $a^{ij}$ is positive definite. Finally we have that 
$$
a^{ij}\varphi_{ik} \varphi_j^k=a^{ij}\sigma^{kl}\varphi_{ik}\varphi_{jl}\geq 0
$$ 
because, as showed in \cite{Di}, if $A$, $B$ and $C$ are symmetric matrices, with $A$ and $B$ positive definite, then $tr(ACBC)\geq 0$.
The thesis follow by the maximum principle. 
\cvd

We can use this partial result to prove that the mean convexity is preserved.

\begin{Lemma}\label{H_bounded}
There exist two positive constants $c_1$ and $c_2$ depending only on $n$ and the initial datum such that for any time the flow is defined
$$
0<c_1\leq H\leq c_2.
$$
\end{Lemma}
\proof From Lemma \ref{eq_evoluz} and the fact that $\left|A\right|^2\geq\frac{1}{4n-1}H^2$ we can compute
$$
\frac{\partial H}{\partial t}\leq \frac{\Delta H}{H^2}-\frac{H}{4n-1}+\frac{4(n+2)}{H}.
$$
By the maximum principle, it is easy to show that $H$ is bounded from above. To prove that H is bounded from below we adapt the proof of Lemma 5.4 of \cite{Pi2} in the new setting: this time we define $\psi=\frac{v}{\sih H}e^{\frac {t}{2(2n+1)}}=\frac{1}{F}e^{\frac {t}{2(2n+1)}}=\frac{\partial\varphi}{\partial t}e^{\frac {t}{2(2n+1)}}$. We want to prove that this function is bounded from above. Preceding like in the proof of Lemma \ref{gradphi_w}:

\begin{eqnarray*}
\frac{\partial \psi}{\partial t} & = & \pat\left(\frac{\partial\varphi}{\partial t}e^{\frac {t}{2(2n+1)}}\right)\\
&=&-\frac{1}{F^2}\left(-a^{ij}\frac{\partial\varphi_{ij}}{\partial t}+a^i\frac{\partial\varphi_i}{\partial t}+\frac{\partial F}{\partial \varphi}\frac{\partial\varphi}{\partial t}\right)e^{\frac {t}{2(2n+1)}}+\frac{1}{2(2n+1)}\psi\\
&=& -\frac{1}{F^2}\left(-a^{ij}\psi_{ij}+a^i\psi_i+\frac{\partial F}{\partial\varphi}\psi\right)+\frac {1}{2(2n+1)}\psi
\end{eqnarray*}

From \eqref{stima01} we have that $\frac{\partial F}{\partial\varphi}\geq 2(2n+1)\frac{\siq}{v^2}$, moreover $\frac{1}{F^2}=\psi^2e^{-\frac {t}{2n+1}}$. By Lemma \ref{gradphi_w} $v^2$ is bounded. Since the function $\siq e^{-\frac {t}{2n+1}}$ is bounded too, we get that

$$
-\frac{1}{F^2}\frac{\partial F}{\partial \varphi}\psi+\frac{1}{2(2n+1)}\psi\leq-c\psi^3+\frac{1}{2(2n+1)}\psi,
$$
for some positive constant $c$. By the maximum principle we deduce that $\psi=\frac{v}{\sih H}e^{\frac{t}{2(2n+1)}}$ is bounded. This imply that there is a constant $c>0$ such that $H\geq cv\frac{e^{\frac{t}{2(2n+1)}}}{\sih}$. The thesis follows since $v\geq 1$ by definition and $\frac{e^{\frac{t}{2(2n+1)}}}{\sih}$ is bounded.
\cvd

Now we can arguing as in Section 5 of \cite{Pi2} and conclude the proof Proposition \ref{gradphi} and prove the asymptotic behaviour of the mean curvature:

\begin{Proposition}\label{H exp}
There is a positive constant $c$ such that
$$
|H-4n-2|^2\leq ce^{-\frac{t}{2n+1}}.
$$
\end{Proposition}

\section{Higher order estimates}\setcounter{equation}{0}
In this section we describe the behaviour of the derivatives of the radial function and its consequence. The results will be only listed. The proofs follow the same procedures of the analogous results of the Sections 6 and 7 of \cite{Pi2}: minor modifications are required when we need to use Proposition \ref{gradphi} and Proposition \ref{H exp}.

\begin{Proposition}\label{higher}
There is a positive constant $c$ such that:
\begin{itemize}
\item[(1)] $|A|^2\leq c$;
\item[(2)] for every $k\in\mathbb{N}$ we have $|\nabla^k_{\sigma}\varphi|^2_{\sigma}\leq ce^{-\frac{t}{2n+1}}$;
\item[(3)] for every $k\in\mathbb{N}$ we have $|\nabla^k_{\sigma}\rho|^2_{\sigma}\leq c$.
\end{itemize}
\end{Proposition}

It follows the uniform parabolicity of equation \eqref{ev_rho} and a uniform $C^2$-estimate for the function $\rho(\cdot, t)$. Arguing as in chapter 2.6 of \cite{Ge1}, we can  apply the $C^{2,\alpha}$ estimates of \cite{Kr} to conclude that the solution of the flow is defined for any positive time and it is smooth, since the initial datum is smooth.

An other interesting consequence of Proposition \ref{higher} is the asymptotic behaviour of the second fundamental form.

\begin{Corollary}
There is a positive constant $c$ such that 
$$
\left|h_i^k-\delta_i^k-\sum_{j=1}^3\delta_i^j\delta_j^k\right|^2\leq c e^{-\frac{t}{2n+1}}.
$$
If we restrict our attention to the horizontal distribution $\mathcal{H}$ (that is we have to consider $i,k\notin\tre$) we get the faster convergence:
$$
\left|h_i^k-\delta_i^k\right|^2\leq c e^{-\frac{2t}{2n+1}}.
$$
\end{Corollary}

An other consequence of Proposition \ref{higher}, part $(3)$ is that, with the same proof of Theorem 8.1 of \cite{Pi2}, we can prove part (3) of Theorem \ref{main}: the rescaled limit metric converges to a $\tilde{g}_{\infty}=e^{2f}\sigma_{sR}$ for some smooth function $f$.

\section{The curvature of the limit metric}

In this section we want to show that the limit metric $\tilde g_{\infty}=e^{2f}\sigma_{sR}$ has not necessarily constant qc-scalar curvature. It is well known that in $\mathbb{HH}^n$ there are no totally umbilical hypersurfaces, so the trace-free part of the second fundamental form cannot have the same strong meaning that it has in the case of hyperbolic space: see  Propositions 3 and 5 of \cite{HW} . Following the ideas developed in \cite{Pi2}, for any star-shaped hypersuface $\mm$ we introduce the following Brawn-York like quantity:
\begin{equation}\label{Q_def}
Q(\mm)=|\mm|^{-1+\frac {1}{2n+1}}\int_{\mm}\left(H-{\hat H}\right)d\mu,
\end{equation}
where $\hat H$ was defined in \eqref{hat_H}. $Q$ gives a measure of how $\mm$ is far from being a geodesic sphere. However $Q$ is not a real measure, because it has not a sign and, even if it is trivially zero when $\mm$ is a geodesic sphere, in general it is not true the contrary. One of the main property of $Q$ is the following.

%Note that we do not know the sign of this function, but, we are sure that $Q(\mm)$ is bounded along the flow 

\begin{Proposition}\label{LQ}
Let $\widetilde{\mm}^{\tau}$ be a family of hypersurfaces in $\hh^n$ that are radial graph of the functions $\tilde{\rho}(z,\tau)=\tau+f(z)+o(1)$, for some fixed $\sss^3$-invariant function $f:\sss^{4n-1}\rightarrow\mathbb R$. Then
\begin{eqnarray*}
\lim_{\tau\rightarrow\infty}Q(\widetilde{\mm}^{\tau})&=&\left(\int_{\sss^{4n-1}}e^{2(2n+1)f}d\sigma\right)^{-1+\frac {1}{2n+1}}\\
&& \phantom{aaa}*\int_{\sss^{4n-1}}e^{(2(2n+1))f}\left(e^{-f}\Delta_{\sigma}e^{-f}-(2n+1)\left|\nabla_{\sigma}e^{-f}\right|^2\right) d\sigma.
\end{eqnarray*}

Moreover if $\lim_{\tau\rightarrow\infty}Q(\widetilde{\mm}^{\tau})\neq 0$, then $e^{2f}\sigma_{sR}$ - the limit of the rescaled metric on $\widetilde{\mm}^{\tau}$ - does not have constant qc-curvature.
\end{Proposition}

The proof is the same, with minor changes concerning the explicit value of the constants, to the proof of Proposition 9.1 of \cite{Pi2}, except for the fact that, in the present case, we have that $\rho_i=0$ for any $i\in\tre$.

If we compare $Q$ with the modified Hawking mass studied for the real hyperbolic case in \cite{HW}, $Q$ has the disadvantage that it works only with $\sss^3$-invariant data and it does not characterize the constant curvature limit. However Proposition \ref{LQ} suggests that, like in the case of complex hyperbolic ambient manifold \cite{Pi2}, the study of the asymptotic behaviour of $Q$ is enough to find a family of initial data such that the limit of the rescaled metric does not have constant qc-curvature. In order to complete this goal we need to study the evolution equation of $Q$ along the inverse mean curvature flow.

\begin{Lemma}\label{eq_evoluz_2}
For any star-shaped $\mm_0$ the following evolution equation holds:
\begin{eqnarray*}
\frac{d Q(\mm_t)}{dt} & = &\frac{1}{2n+1}Q(\mm_t)-|\mm |^{-1+\frac{1}{2n+1}}\int\frac{1}{H}\left(|A|^2-4(n+2)\right)d\mu\\
& & +|\mm |^{-1+\frac{1}{2n+1}}\int\frac{v}{H}\left(\frac{4n-1}{\siq}-\frac{3}{\coq}\right)d\mu.
\end{eqnarray*}
\end{Lemma}
\proof
Since $\hat{H}=(4n-1)\frac{\coh}{\sih}+3\frac{\sih}{\coh}$ and $\frac{\partial \rho}{\partial t}=\frac{v}{H}$, we can easily compute:
$$
\frac{\partial\hat{H}}{\partial t}=\frac{v}{H}\left(\frac{3}{\coq}-\frac{4n-1}{\siq}\right).
$$
The thesis follows using the evolution of $h$ in Lemma \ref{eq_evoluz} and the fact that $$\frac{\Delta H}{H^2}-2\frac{|\nabla H|^2}{H^3}=-\Delta\left(\frac{1}{H}\right),$$ hence its integral vanishes.

\cvd

Now we can prove that if $Q$ decreases, it decreases very slowly.

\begin{Proposition}\label{evoluz_Q}
Let $\mm_t$ an $\sss^3$-invariant star-shaped hypersurface of $\mathbb{HH}^n$ evolving by inverse mean curvature flow. The is a positive constants $c$ such that
\begin{eqnarray*}
\frac{\partial Q}{\partial t} & \geq & -ce^{-\frac{t}{2n+1}}.%-ce^{\frac{t}{2n}}|\nabla H|.
\end{eqnarray*} 
\end{Proposition}
\proof

In case of $\sss^3$-invariant hypersurfaces, by Lemma \ref{hessiani}, \eqref{sff_up} and \eqref{mc} we have
\begin{eqnarray}
\nonumber |A|^2-4(n+2) & = & h_i^kh_k^i-4(n+2)\\
\nonumber& = & \frac{\hat{\varphi}_{ij}\hat{\varphi}_{kh}\tilde{e}^{jk}\tilde{e}^{hi}}{v^2\siq}-\frac{2\coh}{v^2\siq}\hat{\varphi}_{ij}\tilde{e}^{ji}-4(n+2)\\
\nonumber&&(4n-1)\frac{\coq}{v^2\siq}+3\frac{\siq}{v^2\coq}+\frac{6}{v^2}\\
\nonumber & = & \frac{\varphi_{ij}\varphi_{kh}\tilde{\sigma}^{jk}\tilde{\sigma}^{hi}}{v^2\siq}+\frac{6}{v^2}|\nabla_{\sigma}\varphi|^2_{\sigma}\\
&&-\frac{2\coh}{v\sih}\left(H-\hat{H}+\frac{\hat{H}}{v(v+1)}|\nabla_{\sigma}\varphi|^2_{\sigma}\right)\\%\label{d03}
\nonumber& & +\frac{4n-1}{v^2\siq}-\frac{3}{v^2\coq}-\frac{4(n+2)}{v^2}|\nabla_{\sigma}\varphi|^2_{\sigma}.
\end{eqnarray}
Then, by Lemma \ref{eq_evoluz_2}, we have:
\begin{eqnarray*}
|\mm_t|^{1-\frac{1}{2n+1}}\frac{\partial Q(\mm_t)}{\partial t} & = & \int\left(\frac{1}{2n+1}-\frac{2\coh}{vH\sih}\right)\left(H-\hat{H}\right)d\mu\\
&&-\int\frac{\varphi_{ij}\varphi_{kh}\tilde{\sigma}^{ik}\tilde{\sigma}^{jh}}{v^2H\siq}d\mu\\
& & +\int\frac{1}{H}\left(\frac{4n-1}{\siq}-\frac{3}{\coq}\right)\left(v-\frac{1}{v^2}\right)d\mu\\
&&+\int\frac  {|\nabla_{\sigma}\varphi|^2_{\sigma}}{H}\left(4n+2-\frac{2\hat{H}\coh}{(v+1)\sih}\right)d\mu.
\end{eqnarray*}
Since $v$ and $H$ are bounded away from zero, by Proposition \ref{gradphi}, Proposition \ref{H exp} and Proposition \ref{higher} we have that every integral in the equality above is bigger than $-ce^{-\frac{2t}{2n+1}}$ for some positive constant $c$.
\cvd

Now, following the strategy of \cite{HW} and \cite{Pi2} we can complete the proof of Theorem \ref{main}

\begin{Proposition}
There is an $\mm_0$ such that the rescaled induced metric $\tilde{g}_{\infty}$ does not have constant qc-scalar curvature.
\end{Proposition}
\proof
Fix a positive constant $c_0$ big enough and choose an $\sss^3$-invariant function $f:\sss^{4n-1}\rightarrow\mathbb{R}$ such that
$$
\left(\int_{\sss^{4n-1}}e^{2(2n+1)f}d\sigma\right)^{-1+\frac {1}{2n+1}}\int_{\sss^{4n-1}}e^{(2(2n+1))f}\left(e^{-f}\Delta_{\sigma}e^{-f}-(2n+1)\left|\nabla_{\sigma}e^{-f}\right|^2\right) d\sigma%\geq 4c_0,
$$
is bigger than $4c_0$, and consider the family of $\sss^3$-invariant star-shaped hyperfurfaces $\widetilde{\mm}^{\tau}$ defined by the radial function $\tilde{\rho}^{\tau}(z)=\tau+f(z)$. We can fix a $\tau$ big enough such that $\widetilde{\mm}^{\tau}$ is mean convex, and, by Proposition \ref{LQ}, $Q(\widetilde{\mm}^{\tau})\geq 2{c_0}$. Let us consider $\mm_t^{\tau}$, the evolution by inverse mean curvature flow of such $\widetilde{\mm}^{\tau}$. The constant $c$ appearing in Proposition \ref{evoluz_Q} in uniformly bounded if $\tau$ is big enough. Since $c_0$ is big enough too, Proposition \ref{evoluz_Q} ensures that 
$$
\lim_{t\rightarrow\infty} Q(\mm_t^{\tau})\geq {c_0}>0.
$$
The thesis follows from Proposition \ref{LQ}.
\cvd

%\vspace{1cm}

\noindent \textbf{Acknowledgments:} This research is supported by the ERC Avanced Grant 320939, Geometry and Topology of Open Manifolds (GETOM)

\bigskip

%\noindent Giuseppe Pipoli, Institut Fourier, Universit\'e Grenoble Alpes, France.\\ E-mail: giuseppe.pipoli@univ-grenoble-alpes.fr

\end{document}